\documentstyle[12pt]{article}
\begin{document}
%\pagestyle{fancyplain}
%\lhead[\fancyplain{}{\thepage}]{}
%\rhead[]{\fancyplain{}{\thepage}}
%\chead[\fancyplain{}{{\sc Martin Pikaart}}]{\fancyplain{}{{\sc An orbifold
%%%partition of ${\overline{M}_g^n}$}}}
%\cfoot[\fancyplain{\thepage}{}]{\fancyplain{\thepage}{}}

\newtheorem{definition}[subsection]{Definition}
\newtheorem{theorem}[subsection]{Theorem}
\newtheorem{lemma}[subsection]{Lemma}
\newtheorem{proposition}[subsection]{Proposition}
\newtheorem{corollary}[subsection]{Corollary}
\def\square{\Box}
%\newtheorem{example}[subsection]{Example}
%\newremark{examples}[definition]{Examples}
\newtheorem{remark}[subsection]{Remark}
\newtheorem{observation}[subsection]{Observation}
%\newremark{assumption}[definition]{Assumption}
\newcommand{\qed}{\unskip\nolinebreak\hfill\hbox{\quad $\square$}}
\newenvironment{prf}[1]{\trivlist
\item[\hskip
\labelsep{\it #1.\hspace*{.3em}}]}{%~\hspace{\fill}~$\square$%
\endtrivlist}
\newenvironment{proof}{\begin{prf}{Proof}}{\end{prf}}
\hyphenation{para-metrize}

\newcommand\Mgbar { \mbox {${\overline{M}_g}$}}
\newcommand\Mgnbar { \mbox {${\overline{M}_g^n}$}}
\newcommand\ra{\rightarrow}
\hyphenation{stra-ti-fi-ca-tion}

\title{An orbifold partition of \protect{\Mgnbar} }
\author{Martin Pikaart}
\date{}
\maketitle
\begin{abstract} \noindent
We define a partition of $\Mgnbar$ and show that the cohomology of $\Mgnbar$ in
a given degree admits a filtration whose respective quotients are isomorphic
to the shifted cohomology groups of the parts
if $g$ is sufficiently large. This implies that the map $H^k(\Mgnbar) \ra
H^k(M_g^n)$ is onto and that the Hodge structure of $H^k(M_g^n)$
is pure if $g \geq 2k+1$.
The main ingredient is the Stability Theorem of Harer and Ivanov.
\end{abstract}

\section{Introduction}

The moduli space of smooth $n$-pointed complex curves of genus
$g$ is a quasi-projective orbifold $M_g^n$ of dimension $3g-3+n$
(where as usual, we assume that $2g-2+n>0$).
It is compactified by the moduli space of stable pointed curves, $\Mgnbar$,
which is a projective orbifold. We will write $M_g$ and $\Mgbar$ if $n=0$.

We need some terminology in order to state our main result.
Let $X$ be an irreducible orbifold. An {\it orbifold partition} of $X$ is a
finite filtration by closed subvarieties
$$X_\bullet:= (X =X_0 \supset X_1 \supset \cdots \supset X_{m+1} =\emptyset),$$
such that $ Y_\alpha:=X_\alpha \setminus X_{\alpha+1} $ is an orbifold.
If $\dim Y_\alpha >\dim Y_{\alpha +1}$ this is called a
{\it stratification}. A connected component of a $Y_\alpha $ is called a {\it
part} (respectively a {\it stratum}). For example,
a stratification of $\Mgnbar$ is defined by the subvarieties
$X_i:=\{ \mbox{curves with at least $i$ singularities}\}$ for $i=0,1,\dots$;
the strata are just the loci that parametrize stable pointed curves
of a fixed topological type. We will therefore refer to this as the {\it
stratification of $\Mgnbar$ by topological type}. A simple example of an
orbifold partition (of $ {\bf C}^2$) that is not a stratification is given by
${\bf C}^2 \supset L_1 \cup L_2 \supset L_1$, where $L_i$
is the ith coordinate axis; the parts are ${\bf C}^2 \setminus (L_1\cup L_2)$,
$L_2 \setminus \{ 0 \}$ and $L_1$.

Let us say that an orbifold partition $X_\bullet $ of $X$ {\it filters
cohomology up to degree $k$} if the Gysin map
$$H^{j-2 {\it  codim} Y_\alpha}(Y_\alpha)(-{\it codim} Y_\alpha) \rightarrow
H^j(X \setminus X_{\alpha+1})$$
is injective for all $j \leq k$ and all $\alpha$. (Here and throughout this
paper cohomology is taken with rational coefficients.)  For
$2{\it  codim} Y_\alpha > k$ this condition is empty. Thus, if one
is only interested in a partition filtering cohomology up to degree $k$,
it suffices to consider $X_0 \supset \cdots \supset X_{\alpha+1}$ such that
$X_{\alpha+1}$ has codimension at least $\frac{1}{2}k$. It is easily seen that
if a partition filters cohomology up to degree $k$ there
exists a filtration on $H^j(X)$ for $j \leq k$, such that there is a canonical
morphism of mixed Hodge structures which maps the respective
quotients isomorphically onto the groups
$ H^{j -2 {\it codim} Y_\alpha}(Y_\alpha)(-{\it codim} Y_\alpha)
\mbox{ for } j<k \mbox{ and all } \alpha$.
We can now formulate our main result.

\medbreak\noindent{\bf Corollary \ref{coho Mgnbar}}
{\it Given $k \geq 0$, then for $g$ large enough, $\Mgnbar$ admits an orbifold
partition which has $M_g^n$ as its open part, is coarser than the
stratification by topological type and filters cohomology up to degree $k$.}

\medbreak
We shall see that the stratification by topological type does not
have this property. We deduce from the main result the following corollary.

\medbreak\noindent{\bf Corollary \ref{pure Hodge}}
{\it If $g \geq 2k+1$, then the restriction map $H^k(\Mgnbar) \ra H^k(M_g^n)$
is surjective; consequently the mixed Hodge structure on $H^k(M_g)$ is pure of
weight $k$.}

\medbreak
Corollary \ref{coho Mgnbar} enables us to define the ``stable cohomology''
of $\Mgnbar$ as $g$ goes to infinity and $n$ is fixed. This stable cohomology
is not finitely generated. For example, if $n=0$,
the generators in degree 2 are the tautological class and the boundary classes,
naturally indexed by $0,1,2, \dots$.

Using additional properties of the partition that we construct we prove that
for any
$n \geq 0$ the cohomology of $\Mgnbar$ is not of Tate type if $g$ is
sufficiently large.

Let us sketch the proof of corollary \ref{coho Mgnbar}
 and give the motivation for the partition we define. (See
Sections \ref{graphs and partition} and \ref{results} for details.)
We take $n=0$ for simplicity. Let $D_0$ be the locus in $\Mgbar$
parametrizing irreducible curves with one node and $D_1$ the
locus in $\Mgbar$ parametrizing curves with one smooth
component of genus $g-1$ and one of genus 1, joined in
 one point. Consider the maps $f_i:H^{l-2}(D_i) \ra H^l(D_i)$
for $i=0$ or 1, given
by taking the cup product with the first Chern class of the normal bundle
of $D_i$. The normal direction of $D_i$ corresponds to smoothening
the unique singular point. Notice that $f_i$ can be obtained as the composite
of the Gysin map $H^{l-2}(D_i)(-1) \ra H^l(M_g \cup D_i)$ and
the obvious restriction map.
The unique singular point lies in both cases on at least one local component
of high genus.
As we shall see in Section \ref{results}, this will imply that
 the $f_i$ are injective in low degree and
a fortiori the Gysin maps are injective.

However, this does not carry over to codimension 2. Let $D_{01}$
be the locus in $\Mgbar$ parametrizing curves with two nodes
and two components,
one smooth component of genus $g-1$ and one singular component of geometric
genus 0. Clearly one of the nodes lies on two local components of genus 0.
 As we shall see, this implies that $H^{l-4}(D_{01})(-2) \ra H^l(D_{01})$ is
not injective (in fact zero).

A way to restore injectivity is by coarsening the
partition. If we let $D$ be the union of $ D_1$ and $D_{01}$, then the normal
direction of $D$ corresponds to smoothening the singular point on
the component of genus $g-1$ and we get an injective map as we wanted.
Smoothening the singularity on the genus 0 component has become a direction
along the part.

It will turn out that we can define the partition completely in
terms of graphs and certain subgraphs. Our parts will be (irreducible) unions
of strata of the stratification by topological type; locally they are like
the parts in the given example of a non-stratification of ${\bf C}^2$.
The unique open dense stratum of the stratification
by topological type in a given part will be refered to
as its {\it generic} stratum. For every part, its generic stratum
will have a graph with the
property that every edge corresponds to a singularity lying on at
least one local component of high genus.

Unfortunately, I was unable to define such a subgraph canonically; Definition
\ref{specifiek} involves the choice of the number $\alpha$.
The point is that smoothening many singularities lying on low genus components
may yield a component of high genus. See Remark \ref{infinite genus}
for further comments.

 In Section \ref{graphs and partition} we prove the existence of an
orbifold partition using only some formal combinatorial properties of graphs.
In Section \ref{stable and weak} we define weak subgraphs and prove the
combinatorial properties needed in the previous section. Section \ref{results}
contains the results.

\smallskip
I thank Eduard Looijenga for suggesting the problem to me and
for many helpful conversations. I am grateful to the referee for useful
comments.

\section{Stable graphs and the orbifold partition}\label{graphs and partition}

All graphs considered in this paper are weighted, that is,
each vertex has a non-negative integer associated to it.
For a stable $n$-pointed curve $C$, we define its
{\it stable graph} as follows. Each irreducible component of $C$ defines a
vertex, the weight of the
vertex being the genus of the normalisaton of the component.
Omitting from the normalisation of
such an irreducible component the inverse images of the
singular and marked points on $C$ yields a smooth projective curve
minus a finite number of points.
For each missing point we draw a half edge emanating from the vertex defined by
that connected component and for every singular point on $C$ we join the two
half edges associated to it to obtain a whole edge. We omit the word whole if
it is clear that we consider a whole edge.
The remaining half edges, which correspond to the $n$ marked points
on the curve, will be called {\it loose half edges}.
Thus, unless we specify, a half edge can be either
half of a whole edge or a loose half edge.

We define a {\it stable graph}
($n$-pointed of genus $g$) as the stable
graph of some stable $n$-pointed curve of genus $g$.
An automorphism of a stable graph is an automorphism of the underlying
unweighted graph preserving the weights.

Given a stable graph one can recover a topological model for its
stable $n$-pointed curve by taking one curve for every vertex with genus equal
to the weight of the vertex, omitting small open discs for every half edge
emanating from that vertex, glueing the appropriate
boundaries and contracting them. Notice that if we had incorporated the
ordering of the $n$ points in the definition of a stable graph, then our graphs
would correspond bijectively to strata of the stratification by topological
type.

A {\it subgraph} is a subset of a stable graph with the property that if it
contains a half edge it contains its unique vertex.
 By a {\it full subgraph} we mean a
subgraph which contains every half edge emanating from one of its vertices.
A connected component of a full subgraph is a stable graph.
In the same way as above, one can
 associate to a connected component of a full subgraph a topological
model of a curve. We define the {\it genus} of a connected component of a full
subgraph as the arithmetic genus of the curve associated to that connected
component.
The {\it genus} of a full subgraph is defined as the sum of the
genera of its connected components. In other words, the genus of a full
subgraph equals the first Betti number of that full subgrabh plus the sum
of the weights taken over the vertices in the full subgraph.

We denote the greatest integer less than or equal to a real number $a$ by
$[a]$.
Fix $k$ and let ${\cal G}_{k}(g,n)$ be the set of
(isomorphism classes of)
stable graphs, $n$-pointed of genus $g$ and at most $[k/2+1]$ whole edges.
If there is no chance of confusion, we write ${\cal G}$ for
${\cal G}_{k}(g,n)$.
(For the reason we only consider graphs with this many edges,
see the remark made in the introduction
following the definition of a orbifold partition which filters cohomology.)
The maximal number of vertices for a graph
in ${\cal G}$ is $[k/2+1]+1$. From now on we will only consider graphs in
${\cal G}$.

We may regard a stable graph as a one dimensional topological space.
If $e$ is a whole edge, we denote the contraction of $e$ by $\pi_e: \Gamma
\rightarrow \Gamma /e$ and the image of its end vertices by $\pi_e(e)$.
We make $\Gamma/e$ into a weighted graph as follows.
If $p$ is a vertex of $\Gamma$ not incident with $e$, then the weight of
$\pi_e(p)$ is that of $p$.
The weight of $\pi_e(e)$ is defined as the sum of the weights of the two
end vertices if $e$ is not a loop and the weight of its unique end vertex plus
one otherwise. Clearly this makes $\Gamma/e$ into an element of
${\cal G}_k(g,n)$, corresponding to the stable curve obtained by smoothening
the singularity corresponding to the edge $e$.
If $\Delta$ is a full subgraph of $\Gamma$, we denote the quotient of
$\Gamma$ obtained by
contracting all whole edges of $\Delta$ by
$\Gamma/ \Delta$.
The image of a full subgraph under $\pi_e$ is full,
provided that
the contracted edge has either both or neither of its end vertices in that
full subgraph.

The permutation group on $n$ elements, $S_n$, acts on $\Mgnbar$
fixing the locus of points parametrizing curves with at least a fixed number
of singularities. Such loci define the stratification by topological type,
see introduction. Thus $S_n$ acts on the stratification by topological type.
The orbits of the $S_n$-action on the stratification by topological type
of $\Mgnbar $ minus loci of codimension at least $[k/2+1]+1$ are in one to
one correspondence with the elements of ${\cal G}_k(g,n)$.
Giving a partition of the first set is therefore equivalent to giving a
partition of ${\cal G}_k(g,n)$.

This, in turn, corresponds to a function
$$\phi : {\cal G}_k(g,n) \ra {\cal G}_k(g,n),$$
as follows. Recall that we want the parts we seek to be irreducible unions of
strata  of
the stratification by topological type. A graph will be mapped to
the graph of the generic stratum of the part it will belong to. (Compare the
example in the introduction: the graph of $\Delta_{01}$ is mapped to
the graph of $\Delta_1$.) This implies that $\phi$ contracts
a certain subgraph.

So we have to assign to every graph in
${\cal G}_k(g,n)$ a subgraph in such a way that the function $\phi$,
defined as contraction of that subgraph, corresponds to an orbifold
partition.
One important ingredient we need to obtain an orbifold partition is a partial
order on the image of $\phi$,
such that the union of all parts greater than or equal to a given part will be
a Zariski open neighbourhood of that part.
Now we will state the formal properties of the graphs needed to
obtain an orbifold partition- which will turn out to filter cohomology.

\begin{definition}
Suppose we are given, for every $\Gamma$ in ${\cal G}_k(g,n)$, a subgraph
$\Gamma_W$. Define $\phi:{\cal G}_k(g,n) \ra {\cal G}_k(g,n)$ by $\phi(\Gamma)
=\Gamma/\Gamma_W$. Denote the image of $\phi$ by $I$.
For $ \Gamma$ in $I$ define $S_\Gamma^0$ to be the locus in $\Mgnbar$ whose
points correspond to stable $n$-pointed curves with
graph $\Gamma$.
Define $S_\Gamma$ (respectively, $S$) to be the locus in $\Mgnbar$ whose
points correspond to stable $n$-pointed curves with graph
in $\phi^{-1}(\Gamma)$ (respectively,  ${\cal G}$).
\end{definition}

If $C$ is a curve with graph $\Gamma$, then we will call the irreducible
components of $C$ corresponding to  vertices in $\Gamma_W$ {\it weak
components}. Analogously we define {\it strong components}.

\begin{proposition} \label{formal}
Notations as in the previous definition.
Suppose a partial order is given on $I$, such that the following properties
hold:
\begin{enumerate}
\item[a.] $\Gamma_W$ is a proper full subgraph of $\Gamma$ which is invariant
under the automorphism group of the stable graph $\Gamma$;
\item[b.] if $e$ is a whole edge of $\Gamma$ not contained in $\Gamma_W$,
then $\Gamma/e > \Gamma$;
\item[c.] if $e$ is a whole edge contained in $\Gamma_W$ , then the image of
$\Gamma_W$ under contraction of $e$ is $(\Gamma/e)_W$;
\item[d.] a whole edge $e$ is contained in $\Gamma_W$ if and only if
$\pi_e(e)$ is contained in $(\Gamma/e)_W$;
\item[e.] if $k=1$ (respectively $k \geq 2$) and $p$ a vertex of $\Gamma$ such
that the weight of $p$ is 0 (respectively the weight of $p$ is at most
$k+2$), then $\Gamma_W$ contains the vertex $p$.
\end{enumerate}
Then the following hold :
\begin{enumerate}
\item $S_{{\it max }(I)} = M_g^n$;
\item for all $\Gamma \in I$: (the connected components of) $S_\Gamma$ and
$S_\Gamma^0$ are
orbifolds; the generic curve $C_\Gamma$ has graph $\Gamma$;
\item for all $\Gamma \in I$: $\cup_{\Delta \geq \Gamma} S_\Delta$ is a Zariski
open neighbourhood of $S_\Gamma$;
\item the (connected components of the) $S_\Gamma$ are parts of an orbifold
partition;
\item for all $\Gamma \in I$:
the strong (respectively, weak) components of $C_\Gamma$ do not (respectively,
do) degenerate in $S_\Gamma$;
\item for all $\Gamma \in I$:
every singular point of $C_\Gamma$ lies on at least one strong component
of $C_\Gamma$;
\item for all $\Gamma \in I$:
if $C'$ is a strong component of $C_\Gamma$ and $k=1$ (respectively $k \geq 2$)
then the genus of $C'$ is at least 1 (respectively at least $k+3$).
\end{enumerate}
\end{proposition}

\begin{proof}
{\it 1}. Notice that by property {\it a}, $\Gamma_W$ is a proper full subgraph,
which means that $\Gamma/\Gamma_W$ has no whole edges if and only if $\Gamma$
has no whole edges.
By property {\it b}, the stable graph of a smooth $n$-pointed curve
is larger than every other graph; consequently, $S_{max(I)}=M_g^n$.

{\it 2}. We already remarked that (a connected component of) $S_\Gamma^0$ is an
orbifold, namely it is a stratum
of the stratification by topological type. By the definition of $\phi$
we have that $S_\Gamma$ is contained in the closure of $S_\Gamma^0$.
So it suffices to prove that it is locally closed and without self
intersection. It follows from property {\it c}
that $S_\Gamma$ is locally closed.
Suppose $S_\Gamma$ has self intersection; this is necessarily a normal
crossing.
Let $\Delta$ be the graph corresponding to the generic locus of self
intersection of $S_\Gamma$ and suppose for simplicity that there are two local
branches of $S_\Gamma$. Thus, locally in a neighbourhood of $S_\Delta$,
there is an involution permuting the two branches of $S_\Gamma$.
This involution induces an involution $i$ on $\Delta$ which keeps $\Delta_W$
fixed but not pointwise fixed. If $e$ is an edge of $\Delta_W$ such that
$\Delta/e$ is the graph of one of the local branches, then $\Delta/i(e)$
is the (isomorphic) graph of the other local branch. But by property {\it a},
$\Delta_W$ contains both $e$ and $i(e)$, so $S_\Gamma$ contains the plane
spanned by the two local branches, contradiction.

{\it 3}. This follows directly from statement {\it 2} and property {\it b}.

{\it 4}. Choose an order reversing injective map $m: I \ra {\bf N}$ and
define $X_i:= \cup_{m(\Gamma) \leq i} S_\Gamma$.
The $X_i$ form a filtration which defines an orbifold partition by statements
{\it 2} and {\it 3} whose parts are the $S_\Gamma$.

{\it 5}. Follows from property {\it d}.

{\it 6}. By property {\it c}, we have that the image of $\phi$ equals the fix
point set of $\phi$, which means that for $\Gamma$ in the image of $\phi$,
every whole
edge of $\Gamma$ has at least one end vertex not in $\Gamma_W$.

{\it 7}. This is clear from property {\it e}.
\qed \end{proof}

Now we have reduced the problem of finding a partition to a purely
combinatorial one, which will be dealt with in the next section.

\section{Stable graphs, weak subgraphs}\label{stable and weak}

We will now start out to define a full subgraph $\Gamma_W$ for every
$\Gamma$ in ${\cal G}_k(g,n)$ and prove the properties a,b,c,d and e of
Proposition \ref{formal}.

For any full subgraph $\Delta$ of a stable graph $\Gamma$,
we denote its number of
vertices respectively its genus by $v(\Delta)$ resp. $g(\Delta)$. We denote
the weight of a vertex $P$ by $w(P)$.
Let $k \in {\bf Z}_{\geq 0}$ be given.
Let $\Phi_k :{\bf Z}_{\geq 0} \ra {\bf R}_{\geq 0} $ be a function
satisfying
$$2 \Phi_k(n) +[\frac{1}{2}k+1] \leq \Phi_k(n-1), \mbox{ if } n > 0.$$
In particular, it is a decreasing function.

\begin{proposition} \label{union} Let $k \geq 0$ be given.
 Let $\Delta_1$ and $\Delta_2$ be full subgraphs of a stable graph $\Gamma$.
Let $\Delta$ be the full subgraph on the vertices of $\Delta_1$ and
$\Delta_2$.
If $g(\Delta_i) \leq \Phi_k(v(\Gamma)-v(\Delta_i)) $ for $i=1$ and $i=2$,
then also $g(\Delta) \leq \Phi_k(v(\Gamma)-v(\Delta)) $.
\end{proposition}

\begin{proof} We may suppose $ v(\Delta_1) \geq 1$ and $ v(\Delta_1) \geq
v(\Delta_2)$.
The statement is trivial if $\Delta_2$ is contained in $\Delta_1$, so we may
assume that this is not the case. Thus $v(\Delta) > v(\Delta_1)$.
We have:
$$ g(\Delta) \leq g(\Delta_1) + g(\Delta_2) +[\frac{1}{2}k+1] \leq
2 \Phi_k(v(\Gamma)-v(\Delta_1)) +[\frac{1}{2}k+1] $$
$$ \leq \Phi_k(v(\Gamma)-v(\Delta)). \eqno{\Box}$$
\end{proof}

Here and in the next proposition the term $[\frac{1}{2}k+1]$ comes from the
fact that we take a {\it full} subgraph; compare the formula for the genus
of a full subgraph given in Section \ref{graphs and partition}.

\begin{definition} \label{weak}
A full subgraph $\Delta$ of a stable graph $\Gamma$ is called
$\Phi_k$-weak if $g(\Delta) \leq \Phi_k(v(\Gamma)-v(\Delta))$.
\end{definition}

Notice that this definition only depends on $k$ and not on $g$.
By the proposition above the maximal $\Phi_k$-weak subgraph is well defined.

\begin{definition} \label{maximal weak}
Denote the maximal $\Phi_k$-weak subgraph of $\Gamma$ by $\Gamma_W$.
A vertex is called
strong if it is contained in the complement of $\Gamma_W$. Denote by $\Gamma_S$
maximal full subgraph on the strong vertices.
\end{definition}

Not every full subgraph of the maximal $\Phi_k$-weak subgraph is $\Phi_k$-weak,
see the example following Proposition \ref{abcde hold}. Therefore we don't use
the term $\Phi_k$-weak vertex.

\begin{proposition}\label{full proper invariant}
If $g > \Phi_k(0)$, then for every $\Gamma$ in ${\cal G}_k(g,n)$ we have that
$\Gamma_W$ is a proper full subgraph which is invariant under
the automorphisms of the stable graph $\Gamma$.
\end{proposition}

\begin{proof}
By definition, $\Gamma_W$ is full. It is a proper subgraph because
$g(\Gamma)>\Phi_k(v(\Gamma)-v(\Gamma))$ and thus $\Gamma$ is not
$\Phi_k$-weak. Invariance under automorphisms follows because both
the genus and the number of vertices are kept fixed by automorphisms
of the stable graph.
\qed \end{proof}

Recall that if $e$ is a whole edge of $\Gamma$, we write $\pi_e$ for the
contraction of $e$ and $\pi_e(e)$ for the image of its end vertices.
The inverse image of a full subgraph under $\pi_e$ is full.
For $e$ a whole edge of $\Gamma$, we put $\delta(e):=0$ if $e$ is a
loop, $\delta(e):=1$ otherwise.

\begin{proposition} \label{sub}
If $e$ is a whole edge contained in $\Gamma_W$ or in its complement,
then we have $\pi_e(\Gamma_W) \subset (\Gamma /e)_W$.
\end{proposition}

\begin{proof} We have:
$$g(\pi_e(\Gamma_W))=g (\Gamma_W) \leq \Phi_k(v(\Gamma)-v(\Gamma_W)) =$$
\[ \left\{ \begin{array}{ll}
  \Phi_k((v(\Gamma/e) +\delta(e)-v(\pi_e(\Gamma_W))) &
\mbox{  $e$ outside $\Gamma_W$}\\
\Phi_k(v(\Gamma/e) +\delta(e)-[v(\pi_e(\Gamma_W))+\delta(e)]) &
\mbox{  $e$ contained in $\Gamma_W$}
\end{array} \right. \]
$$ \leq \Phi_k(v(\Gamma/e)-v(\pi_e(\Gamma_W))). \eqno{\Box}$$
\end{proof}

\begin{proposition} \label{gammaw maps to gamma/ew}
If $e$ is a whole edge contained in $\Gamma_W$, then the image of
$\Gamma_W$ under contraction of $e$ is $ (\Gamma /e)_W$.
\end{proposition}

\begin{proof}
Proposition \ref{sub} yields us one inclusion.
Put $\Delta :=\pi_e^{-1}((\Gamma/e)_W)$, then we have:
$$g(\Delta) =g ((\Gamma/e)_W) \leq \Phi_k(v(\Gamma/e)-v((\Gamma/e)_W)) $$
 $$=\Phi_k(v(\Gamma)-\delta(e)-[v(\Delta)-\delta(e)])
=\Phi_k(v(\Gamma)-v(\Delta)). \eqno{\Box}$$\end{proof}

\begin{proposition}\label{not weak kept fixed under degeneration}
 A whole edge $e$ is contained in $\Gamma_W$ if and only if
$\pi_e(e)$ is contained in $(\Gamma/e)_W$.
\end{proposition}

\begin{proof}
The implication "$\Rightarrow$" follows from Proposition \ref{gammaw maps to
gamma/ew}. For
the other implication, suppose $(\Gamma/e)_W$ contains $\pi_e(e)$.
Put $\Delta :=\pi_e^{-1}((\Gamma/e)_W)$, then the same computation as in the
proof of Proposition \ref{gammaw maps to gamma/ew} implies that $\Delta$ is
 $\Phi_k$-weak and therefore
contained in $\Gamma_W$. Consequently, $\Gamma_W$ contains $e$.
 \qed \end{proof}

\begin{definition} \label{phi}
Define a function $\phi_k(g,n): {\cal G}_k(g,n) \rightarrow {\cal G}_k(g,n)$
 by $\phi_k(\Gamma) =\Gamma /\Gamma_W$. When there is no chance of
confusion, we will write $\phi$ instead of $\phi_k$.
\end{definition}

We call $\phi$ contraction of the maximal $\Phi_k$-weak subgraph.
It is clear that we have: $\Gamma$ is contained in $Im(\phi)$
if and only if $ \Gamma_W$
does not contain any whole edges. It follows that the image of $\phi$ equals
the fix point set of  $\phi$.
Put $I_k(g,n):=Im(\phi_k(g,n))$; we will write $I$ instead of $I_k(g,n)$
if there is no chance of confusion. The set $I$ will be the index set for our
orbifold partition.
Before we can define a partial order on $I$, we need two more propositions.

\begin{proposition} \label{low inside}
If $P$ is a vertex of $\Gamma$ and $w(P) \leq {\mbox max}
 \{ w(Q) \}$, where $Q$ runs over the vertices of $\Gamma_W$, then $\Gamma_W$
contains $P$.
\end{proposition}

\begin{proof}
Suppose not, let $\Delta$ be the full subgraph on $P$ and $\Gamma_W$.
We have :
$$ g(\Delta) \leq g(\Gamma_W)+g(P) +[\frac{1}{2}k+1] \leq 2g(\Gamma_W) +
[\frac{1}{2}k+1] $$
$$ \leq \Phi_k(v(\Gamma)-(v(\Gamma_W) +1))
=\Phi_k(v(\Gamma)-v(\Delta)). \eqno{\Box}$$
\end{proof}

\begin{proposition}\label{joins}
Let $e$ be an edge which joins $\Gamma_S$ and $\Gamma_W$.
Let $\Lambda $ be the full subgraph of $\Gamma/e$ on the vertices of
$\pi_e(\Gamma_W) -\pi_e(e)$. Then $\Lambda$ is $\Phi_k$-weak.
\end{proposition}

\begin{proof} We have:
$$ g(\Lambda) \leq g(\Gamma_W) \leq \Phi_k(v(\Gamma)-v(\Gamma_W))
=\Phi_k(v(\Gamma)-[v(\Lambda)+1])$$
$$=\Phi_k(v(\Gamma/e)-v(\Lambda)). \eqno{\Box}$$
 \end{proof}

\begin{definition} Let $\Gamma$ be in ${\cal G}$.
Let $s(\Gamma)$ be the number of strong vertices.
Define $\{ w_i(\Gamma) \}_{i=1}^{s(\Gamma)}$ to be the set of weights of
strong vertices, ordered such that $w_i(\Gamma) \geq w_{i+1}(\Gamma)$.
For $i  \in \{ s(\Gamma)+1,...,[k/2+1]+1 \}$, define $w_i(\Gamma):=g$.
Let $n(\Gamma)$ be the number of half edges parting from strong vertices.
Finally, define {\rm index}$(\Gamma)$ of a graph $\Gamma$ to be the vector
$(w_1(\Gamma),....,w_{[k/2+1]+1}(\Gamma),n(\Gamma)).$
\end{definition}

Notation: ${\rm index}(\Gamma)> {\rm index}(\Delta)$ refers to lexicographical
ordering.

\begin{proposition} \label{indeks}
a) If $e$ is a whole edge of $\Gamma$ not contained in $\Gamma_W$, then
${\rm index}(\Gamma/e)>{\rm index}(\Gamma)$. \hfill \break
b) If $e$ is a whole edge contained in $\Gamma_W$, then ${\rm index}(\Gamma/e)
={\rm index}(\Gamma)$.
\end{proposition}

\begin{proof}
a)  Suppose first that neither of the end vertices of $e$ is in $\Gamma_W$
and suppose that ${\rm index}(\Gamma)  $ $=  (\dots,w_i , \dots , w_j,\dots)$.
 If $e$ is a loop and its end vertex has weight $w_i$,
then by Propositions \ref{low inside} and \ref{not weak kept fixed under
degeneration} we have
${\rm index}(\Gamma/e)=(\dots,w_i+1,\dots) > {\rm index}(\Gamma)$.
If $e$ is not a loop and its end vertices have weights $w_i$ and $w_j$,
 then ${\rm index}(\Gamma/e)=(\dots,w_i+w_j,\dots) > {\rm index}(\Gamma)$.

Secondly suppose $e$ has end vertices $P \in \Gamma_S$ and $Q \in \Gamma_W$.
By Proposition \ref{joins} we have $s(\Gamma) \geq s(\Gamma/e)
$. If $s(\Gamma) > s(\Gamma/e)$, then we are done.
If $s(\Gamma) = s(\Gamma/e)$ and $w(Q) >0$ the argument above applies.
If $s(\Gamma) = s(\Gamma/e)$ and $w(Q)=0$, then the first $[k/2+1]+1$
coefficients of both indices are equal.
We define the following numbers:
\[ \begin{array}{l}
a:=\# \{ \mbox{edges joining $P$ and $Q$} \},\\
b:=\# \{ \mbox{loops at $Q$} \},\\
c:=\# \{ \mbox{edges joining $Q$ and another  vertex in $\Gamma_W$} \},\\
d:=\# \{ \mbox{edges joining $Q$ and a vertex in $\Gamma_S$ not equal to
$P$}\}, \\
f:=\# \{ \mbox{loose half edges at $Q$}\}.
\end{array} \]
One sees easily that $n(\Gamma/e)=n(\Gamma) +a+2b+c+d+f-2$.
Stability of the vertex $Q$ implies $a+2b+c+d+f\geq 3$, so we have
${\rm index}(\Gamma/e)> {\rm index}(\Gamma)$. \hfill \break
 b) Follows immediately from Proposition \ref{gammaw maps to gamma/ew}.
\qed \end{proof}

\begin{corollary}
The index is preserved under contraction of the maximal $\Phi_k$-weak
subgraph.
\end{corollary}

\begin{definition}\label{partial order}
We define a partial order on $I$ as follows:
$\Gamma \geq \Delta$ if and only if ${\rm index}(\Gamma)$ $ >
{\rm index}(\Delta)$ or $\Gamma =\Delta$.
\end{definition}

\begin{corollary} \label{order OK}
If $e$ is a whole edge of $\Gamma$ not contained in $\Gamma_W$, then $\Gamma/e
> \Gamma$.
\end{corollary}

Now we are almost ready to prove the remaining property e of Proposition
\ref{formal}. We do this by making a suitable choice for $\Phi_k$.
For a given $k$, define $ \alpha := 1/ [k/2+1]+2$ and
$\beta:=[k/2+1]^2 \alpha^{[k/2+1]+1}$. We fix $\alpha$ and $\beta$
for the rest of this paper.

\begin{definition} \label{specifiek}
Define $\Phi_k :{\bf Z}_{\geq 0} \ra {\bf R}_{\geq 0} $ by $\Phi_k(n):=
\alpha^{-n}\beta$.
\end{definition}

An easy calculation yiels that
$$2 \Phi_k(n) +[\frac{1}{2}k+1] \leq \Phi_k(n-1), \mbox{ if } n > 0.$$

\begin{proposition}\label{gammaW bevat laag geslacht}
If $k=1$ (respectively $k \geq 2$) and $P$ a vertex of $\Gamma$ such
that the weight of $P$ is 0 (respectively the weight of $P$ is at most
$k+2$), then $\Gamma_W$ contains the vertex $P$.
\end{proposition}
\begin{proof}
Let $l$ be the number of loops at the vertex $P$. One has to check that
$0+l \leq  \Phi_k([k/2+1]-l) $ (respectively $k+2+l \leq \phi_k([k/2+1]-l)$),
which is an easy computation.
\qed \end{proof}

\begin{proposition} \label{abcde hold}
Notations as above.
Let $k$ be given, and let $g$ be at least $\beta$. Define $\Phi_k$ as in
Definition \ref{specifiek}. Define $\phi$ and $I$ as in Definition
\ref{phi}. Let a partial order on $I$ be defined by Definition \ref{partial
order}. Then the properties a, b, c, d and e of Proposition \ref{formal}
hold.
\end{proposition}

\begin{proof}
The properties a,  b, c, d and e are precisely the Propositions
\ref{full proper invariant}, \ref{order OK}, \ref{gammaw maps to gamma/ew},
\ref{not weak kept fixed under degeneration} and \ref{gammaW bevat laag
geslacht}. \qed
\end{proof}

Not every full subgraph of the maximal $\Phi_k$-weak subgraph is
$\Phi_k$-weak.
Consider the
following example: $k=9$, so $[k/2+1]=5$ and $\alpha=11/5$.
We have $\beta=11^6/5^4 \approx 2834$ and thus we have to take $g > 2834$.
Let $\Gamma$ be the graph with vertices $p_1, \dots , p_6$ and and
edges $e_1, \dots , e_5$ such that $e_i$ joins $p_i$ and $p_{i+1}$.
Suppose $w(p_1)=w(p_2)=11,~w(p_3)=22,~w(p_4)=44,~w(p_5)=88$ and
$w(p_6)=g-176$. One can check that the maximal $\Phi_k$-weak subgraph is the
full subgraph on the vertices $p_1, \dots,p_5$, but the full subgraph on
the vertex $p_5$ is not $\Phi_k$-weak.

\begin{remark} \label{infinite genus} {\rm
If one is willing to consider graphs of infinite genus, i.e.\ several
vertices can have infinite genus, then a canonical
definition of weak subgraph is readily available: just take the full subgraph
on the vertices of finite genus. (The definition of index
needs to be adapted.) Properties analogous to those of the previous
propositions, and in some cases even stronger results, can be proved.}
\end{remark}

\section{The results}\label{results}

We keep the notations and assumptions of Proposition \ref{abcde hold}, unless
the contrary is explicitely stated.

As we have seen in Section \ref{graphs and partition}, a connected component of
the $S_\Gamma^0$ is a stratum of the
stratification by topological type; it is an orbifold of
codimension equal to the number of singular points of its topological model
which we will denote by $C_\Gamma$.
In fact, a  stratum is isomorphic to a product of lower dimensional moduli
spaces $M_h^m$ modulo a finite group.
Recall that the symmetric group on $n$ elements acts transitively on the
connected components of the $S_\Gamma^0$.
 Furthermore, the orbifold
normal bundle of $S_\Gamma^0$ has fibre over $(C,x_1, \dots, x_n)$
isomorphic to $\oplus(T_xC' \otimes T_xC'')$, where the sums runs over
the singular points of $C$ and $C',~C''$ are the local components of $C$ in a
suitable neighbourhood  of $x$. This isomorphism globalizes to an
isomorphism of bundles on $S_\Gamma^0$.

It follows from Proposition \ref{formal} that the orbifold
normal bundle of $S_\Gamma$ has fibre over $(C,x_1, \dots, x_n)$ isomorphic to
$\oplus(T_xC' \otimes T_xC'')$. Here the sums runs over the singular points
of $C$ which are specializations of singular points on the topological
model of $S_\Gamma^0$ and $C',~C''$ are the local components of $C$ in
a suitable neighbourhood  of $x$. Thus, the sums runs over the edges of
$\Gamma$.

Notice that $S=\cup_{\Gamma \in I} S_\Gamma$ and that the real codimension
of the complement of $S$ in $\Mgnbar$ is larger than $k$.

Before going on we explain the Stability Theorem of Harer and Ivanov in an
algebro-geometric way. These theorems
are essential in the proof of Theorem \ref{filters coarse}.

Let $S$ be the locus in $\overline{M_{g+1}}$ of irreducible stable curves
of genus $g$ with one singularity; $S$ is a stratum of the stratification by
topological type, it is isomorphic to  $M_g^2$ modulo the involution permuting
the two points. There is a map $p:S \rightarrow M_g $
which forgets the two points.
Let $i:U_S \hookrightarrow \overline{M_{g+1}}$ be the inclusion of
a suitable $C^\infty$ tubular neighbourhood of $S$ and let $\pi :U_S
\rightarrow S$ be
the natural retraction. Consider the diagram
$$M_{g+1} \stackrel{i}{\longleftarrow} U_S \setminus S \stackrel{p \circ
\pi}{\longrightarrow} M_g.$$
The Stability Theorem of Harer and Ivanov now says that if $g \geq 2k+1$, then
 both maps induce isomorphisms on cohomology in degree up to $k$ (see
\cite{Harer1},\cite{Ivanov}). Furthermore $H^0(M_g) \cong {\bf Q}$  for $g \geq
0$ and $H^1(M_g) \cong 0$ for $g \geq 1$, see \cite[Ch. 7]{Harer2}.
These facts account for the conditions of property e of Proposition
\ref{formal}.

Hence we can define the kth stable cohomology group $H^k(M_\infty)$ of the
moduli space by $H^k(M_\infty):=H^k(M_g)$ when $g \geq 2k+1$.
If $g \geq 2k+1$ (respectively $g \geq 0$ if $k=0$, respectively $g>0$ if
$k=1$), we say $g$ is in the stable range with respect to $k$.
Moreover $i^*$ and $(p \circ \pi)^*$ are morphisms of mixed Hodge structures,
so $H^k(M_\infty)$ carries a well-defined mixed Hodge structure.

\begin{remark} {\rm
In \cite{Mumford} classes $\kappa_i \in H^{i,i}(M_\infty)$ are constructed
and in \cite{Miller} it is proved that the symmetric algebra on these classes
injects into $H^\bullet(M_\infty)$. Mumford conjectures that this is
actually an isomorphism in low degree, see \cite[Introduction]{Mumford}.
The conjecture would imply our corollary \ref{pure Hodge}. } \end{remark}

We need a corollary of Harer's results.
Assigning to a pointed curve the tangent space at its ith point defines
a line bundle ${\cal L}_i$ on $\Mgnbar$.
Consider the natural forgetful map
$M_g^n \rightarrow M_g$. Define $H^\bullet (M_g)[u_1, \dots, u_n]
\rightarrow H^\bullet(M_g^n)$, where the $u_i$ have degree 2,
by sending $u_i$ to the first Chern class
$c_1({\cal L}_i)$. Then we have (see \cite{Looijenga}): this is an isomorphism
up to degree $k$ if $g \geq 2k+1$.

\begin{theorem} \label{filters coarse}
Let $k$ be given. If $g > \beta ,$ then the partition $S=\cup_{\Gamma \in I}
 S_\Gamma$ is coarser than the stratification by topological type,
has $S_{{\it max }(I)} =M_g^n$ and filters cohomology up to degree $k$.
\end{theorem}

\begin{proof}
We have already seen in Proposition \ref{formal}, properties 1 and 4,
that the $S_\Gamma$ form a partition which has $M_g^n$ as open
part and is coarser than the stratification by topological type.
So it remains to show that for all $\Gamma \in I$ and all $l<k$, the Gysin maps
on cohomology $H^{l-2codim S_\Gamma}(S_\Gamma)(-codim S_\Gamma)
\ra H^l(\cup_{\Delta \geq \Gamma} S_\Delta)$
induced by the inclusions $S_\Gamma \rightarrow
\cup_{\Delta \geq \Gamma} S_\Gamma$ are injections.
When we write property x we mean property x of Proposition \ref{formal}.

The orbifold $S_\Gamma^0$ is the quotient of
$$ \prod M_{g_s}^{n_s} \times \prod M_{g_w}^{n_w}$$
by a finite group. Here the first product runs over the strong vertices
of $\Gamma$ and the second over those which are not strong.
 Because of properties 2 and 5 we have that $S_\Gamma$ is contained in the
quotient of
$$\prod M_{g_s}^{n_s} \times \prod \overline{M_{g_w}^{n_w}}$$
 by a finite group, where again the products
run over the vertices which are strong respectively not strong.

Property 7 implies that the $g_a$ are in the stable range w.r.t. $[k/2+1]$.
By the result of Looijenga mentioned above we get :
$$H^i(M_{g_s}^{n_s}) \cong \mbox{degree $i$ part of }H^i(M_\infty)[u_1, \dots,
u_{n_s}], \eqno{(1)}$$
where the $u_i$ have degree 2.

We have seen that the normal bundle of $S_\Gamma$ splits as a direct
sum of line bundles, and thus its top Chern class becomes the product of the
first Chern classes of these line bundles.
We claim that these first Chern classes are all of the form $u_i+u_j$ or
$u_i+a$, where the $u_i$ are as in $(1)$ and $a$ is
an element of $H^*(\overline{M_{g_w}^{n_w}})$. We postpone the proof of the
claim for a moment.

By properties 2 and 3, $S_\Gamma$ is a closed suborbifold of
$\cup_{\Delta \geq \Gamma}S_\Delta$.
Consider the Gysin sequence for the inclusion $S_\Gamma$ in $\cup_{\Delta \geq
\Gamma}S_\Delta$:
$$
 \dots \ra H^{l-2 {\it codim} S_\Gamma}(S_\Gamma)(-codim S_\Gamma)
 \ra H^l(\cup_{\Delta \geq \Gamma} S_\Delta) \ra
H^l(\cup_{\Delta > \Gamma} S_\Gamma) \ra \dots
$$
Composing $H^{l-2codimS_\Gamma}(S_\Gamma)(-codim S_\Gamma) \ra H^l(\cup_{\Delta
\geq \Gamma} S_\Delta)$
with the restriction morphism to $  H^l(S_\Gamma)$
we get a morphism $H^{l-2codimS_\Gamma}(S_\Gamma)(-codim S_\Gamma)
 \ra H^l(S_\Gamma)$
which is given by taking the cup product with the top Chern class of
the normal bundle of $S_\Gamma$. From what we have said above it follows that
cupping with this Chern class is injective up to degree $k$. A fortiori
the Gysin maps are injective.

It remains to prove the claim.
As explained, the line bundles of which we are taking the first Chern classes
correspond to whole edges of the graph $\Gamma$. Let $f$ be a whole
edge and let $P$ and $Q$ be its end vertices (which possibly coincide).
The line bundle under consideration is the tensor product of the two
line bundles corresponding to $P$ and half of $f$ respectively to $Q$ and
the other half of $f$.
By property 6 we have that either $P$ and $Q$ are both strong
vertices or precisely one of them is a strong vertex. In the first case the
first Chern classes of both line bundles are of the form $u_i$ and so
we get $u_i+u_j$ as first Chern class of the tensor product. In the second case
we get $u_i+a$ for some element $a$ of $H^*(\overline{M_{g_w}^{n_w}})$.
This proves the claim.
\qed \end{proof}

\begin{remark} {\rm
The parts $S_\Gamma$ depend upon
the definition of $\Phi_k$-weak subgraphs, which depends upon $\alpha$, which
in turn depends upon $k$.
This implies that the filtration in Theorem \ref{filters coarse}
 depends upon $k$. In this remark we will write $\alpha_k$ and $\beta_k$ to
stress dependence.

There is a natural inclusion $i: {\cal G}_{k}(g,n) \hookrightarrow
{\cal G}_{k+1}(g,n)$ which is a bijection if $k$ is even.
We claim that $i(\Gamma_W) \subset (i(\Gamma))_W$.
This is clear if $k$ is even because then
$\alpha_k^{v(\Delta)-v(\Gamma)}\beta_k
=\alpha_{k+1}^{v(i(\Delta))-v(i(\Gamma))}\beta_{k+1}$ and
we even have equality. If $k$ is odd
one has to check that $\alpha_k^{v(\Delta)-v(\Gamma)}\beta_k
\leq \alpha_{k+1}^{v(i(\Delta))-v(i(\Gamma))}\beta_{k+1}.$
If  we put $l:=[k/2+1]$
and $n=v(\Delta)-v(\Gamma)$,
then this amounts to checking the inequality
 $(\frac{2l+3}{l+1})^{l+1-n}(l+1)^2 \geq (\frac{2l+1}{l})^{l-n}l^2$,
 which is tedious but elementary.
It follows that the parts $S'_\Delta$ for $k+1$ contain unions of parts
$S_\Gamma$ for $k$:
$S'_\Delta \supset \cup_{\Gamma \in J} S_\Gamma$. We can now apply Theorem
\ref{filters coarse} to this union to get a filtration on $H^l(S_\Delta')$ some
of whose
subquotients are isomorphic to $H^{l+2 codimS_\Delta'-2
codimS_\Gamma}(S_\Gamma)(-codim S_\Delta' + codim S_\Gamma)$
for $l<k$. We conclude that enlarging $k$ amounts to taking
the union of some parts.
 } \end{remark}

\begin{remark} {\rm
Theorem \ref{filters coarse} does not hold if we replace the partition by the
stratification by topological type. To see this,
let $k$ be given and take $g>5\alpha^{[k/2+1]+1}$. Take $n=0$ for simplicity.
Choose natural numbers $d>c>b>a>\alpha^{[k/2+1]+1}$ such that
$a+b+c+d=g$. Consider the graph $\Gamma$
which has five vertices, of weights $0,~a~,b,~c,~d$ and four edges,
joining the weight $0$ component to the other four. The automorphism group
of the stable graph $\Gamma$ clearly is trivial. One checks easily that
$\Gamma_W$ is
the full subgraph on the weight $0$ vertex, which implies that $\Gamma$ defines
the open stratum $S_\Gamma^0$ of a part $S_\Gamma$. The part $S_\Gamma$ is
obtained by letting
the genus $0$ curve degenerate in all possible ways. There are three
possible degenerations, corresponding to the three possibilities of
partioning the four other vertices into
two sets of two. So we have: $S_\Gamma^0 \cong N \times M_0^4$ and $S_\Gamma
\cong N \times \overline{M_0^4}$, where $N:=M_a^1 \times M_b^1 \times
M_c^1 \times M_d^1$.
$M_0^4$ is a ${\bf P}^1$ minus 3 points and $\overline{M_0^4}$ is ${\bf P}^1$.
Thus $S_\Gamma=S_\Gamma^0 \cup \cup_1^3 N \times {\it point}$.
 We have $H^2(N \times \overline{M_0^4})=H^2(N) \otimes H^0(\overline{M_0^4})
\oplus H^0(N) \otimes H^2(\overline{M_0^4})=H^2(N) \oplus
{\bf Q}$. Suppose Theorem \ref{filters coarse} would hold with the
partition replaced by the
stratification by topological type. Then we would have, using the above:
$H^2(N \times \overline{M_0^4})=
H^2(N \times {M_0^4}) \oplus \oplus_{i=1}^3 H^0(N \times {\it point})=
 H^2(N) \oplus {\bf Q}^3$, which leads to a
contradiction. }
\end{remark}

\begin{corollary} \label{coho Mgnbar}
 Given $k \geq 0$, then for $g>\beta$, $\Mgnbar$ admits an orbifold partition
which has $M_g^n$ as its open part, is coarser than the stratification by
topological type and filters cohomology up to degree $k$.
\end{corollary}

\begin{proof}
Using Theorem \ref{filters coarse}
 we get the statement with $\Mgnbar$ replaced by $S$.
The real codimension of the complement of $S$ in $\Mgnbar$ is
$2[k/2+1]+2 >k$, thus
$H^l(\Mgnbar) \cong H^l(S)$, for all $l<k$. (Compare the remark following
the definition of a partition which filters cohomology in the introduction.)
Combining these gives the desired result.
\qed \end{proof}

\begin{corollary} \label{pure Hodge}
If $g \geq 2k+1$, then $H^k(\Mgnbar) \ra H^k(M_g^n)$ is onto; consequently,
the mixed Hodge structure on $H^k(M_g^n)$ is pure of weight $k$.
\end{corollary}

\begin{proof}
By corollary \ref{coho Mgnbar}, we have that $H^k(M_g^n)$ is a quotient of
 $H^k(\Mgnbar)$, so the last statement holds for $g > \beta$.
Now use that the image of $H^k(\Mgnbar)$
in $H^k(M_g^n)$ is $W^kH^k(M_g^n)$, which equals $H^k(M_g^n)$ if $g \geq
2k+1$.
\qed \end{proof}

\begin{corollary}\label{not Tate}
If $g >>0$, then the cohomology of $\Mgnbar$ is not of Tate type.
\end{corollary}

\begin{proof}
We use the fact that the modular form $\Delta$ can be seen as a non vanishing
holomorphic 11-form on $\overline{M_1^{11}}$ (see \cite{Deligne}). Thus,
$H^{11,0}(\overline{M_1^{11}})$ is not zero.
Take $g$ large enough and $n$ arbitrary. One of the strata in
corollary \ref{coho Mgnbar}
 is the following: the generic graph consists of two vertices, one of weight
$g-11$ and one of weight 1. There are 11 edges joining them and $n$ loose half
edges at the vertex of weight $g-11$. Clearly, the maximal $\Phi_k$-weak
subgraph of
this graph consists of the vertex of weight 1 with its 11 half edges.
This graph is therefore in the image of
the map $\phi$ and defines a generic point of a part of codimension 11.
The open dense topological stratum in it is $M_{g-11}^{n+11}
\times {M_1^{11}}$.
The corresponding part is obtained by
letting the genus 1 curve degenerate, thus it is $M_{g-11}^{n+11} \times
\overline{M_1^{11}}$.

Consider the cohomology group $H^{22,11}(\Mgnbar)
\subset H^{33}(\Mgnbar)$. By corollary \ref{coho Mgnbar}
 it has a subquotient $H^{11,0}(M_{g-11}^{n+11} \times
\overline{M_1^{11}})$, which has as direct summand
$H^0(M_{g-11}^{n+11}) \bigotimes H^{11,0}(\overline{M_1^{11}}) \neq 0 $.

By considering different strata we find that other cohomology groups don't
vanish either. For example: suppose $n \geq 10$ and take a graph consisting of
two vertices, one
of weight 1 and one of weight $g-1$, joined by one edge. Furthermore the weight
 one vertex has 10 loose half edges and the other vertex $n-10$. This
defines the generic graph of a part of codimension 1.
By the same argument as above, we see that $H^{12,1}(\Mgnbar)$ is not zero.
\qed \end{proof}

\begin{proposition}\label{stabiele cohomologie voor Moneindignstreep}
There exists a ``stable cohomology" of $\Mgnbar$ for $n$ fixed and $g \ra
\infty$.
\end{proposition}

\begin{proof}
Notation as before.
Let $k$ be given and suppose $g > \beta$. Let $\Gamma$ be in the image of
$\phi_k(g,n)$, that is, $\Gamma$ defines a generic stratum of a part $S$.
Choose a strong vertex $P$ of $\Gamma$, let $h$ be the weight of $P$.
Consider the graph we get by changing $h$ to $h+1$. Call this graph
$\Delta$. Because the definition of $\Phi_k$-weak subgraph does not depend on
$g$,
we see that $\Gamma_W$ and $\Delta_W$ can be identified. This means that
$\Delta$ is in the image of $\phi_k(g+1,n)$ and therefore defines a part $T$ in
$\overline{M_{g+1}^n}$. Furthermore we get that $H^l(S) \cong H^l(T)$ for all
$l <k$, because $h$ is in the stable range with respect to $k$.

For every $x \in {\bf N}$
 we define $\psi_x: {\cal G}_k(g+x,n) \ra
{\cal G}_k(g+x+([k/2+1]+1)!,n)$ as follows:
for any vertex $P$ in $\Gamma_S$ add
$([k/2+1]+1)!/ \# \Gamma_S$ to its weight.
$\psi$ clearly respects the parts and
thus defines a map $\psi_x :I_k(g+x,n) \ra I_k(g+x+([k/2+1]+1)!,n)$. Because
the maximal $\Phi_k$-weak subgraphs don't change, this map is an injection.
Furthermore it induces isomorphisms on the cohomology of the parts, as
explained above. So we get an injection induced by $ \psi_x$:
$$H^l(M_{g+x}^n)  \hookrightarrow H^l(M_{g+x+([\frac{1}{2}k+1]+1)!}^n),$$
 for all $l<k$, which maps the subquotients isomorphically onto the
corresponding subquotients. We claim that the inductive limit
over these maps is independent of $x$. We will compare the inductive limits
for 0 and $x>0$. For any $\Gamma$ in $Im(\phi_k(g,n))$ we choose
one strong vertex $P_\Gamma$.
Define $\chi_x(\Gamma) \in {\cal G}_k(g+x,n)$
by adding $x$ to the weight of $P$. $\chi_{x}(\Gamma)$
is in the image of $\phi_k(g+x,n)$. Define $\overline{\chi_x}:
Im(\phi_k(g+([k/2+1]+1)!,n)) \ra
Im(\phi_k(g+x+([k/2+1]+1)!,n))$ by adding $x$ to the
weight of the strong vertex $\psi_x (P_\Gamma)$. We clearly have
$\overline{\chi_x} \circ \psi_0=\psi_x \circ \chi_x$.
Now we replace 0 by $x$ and $x$ by
$([k/2+1]+1)!$.
Playing the same trick and using that the inductive limits are
canonically isomorphic if the indices differ by a multiple of
 $([k/2+1]+1)!$, we prove our claim. We define $H^l(\overline{M^n_\infty})$ as
this inductive limit.
 \qed \end{proof}

One may regard this as the ``stable cohomology" of $\Mgnbar$ for $g \ra
\infty$. We note however that this stable cohomology is not finitely
generated.
For example, if $n=0$, the stable generators in degree $2$ are the tautological
class, and the boundary classes naturally indexed by $0,1,2,\dots$.
We note also that the part defined in Corollary \ref{not Tate} defines
a part for every $g$ sufficiently large; consequently the proof
of Theorem \ref{stabiele cohomologie voor Moneindignstreep} shows that
even the stable
cohomology of $\Mgnbar$ for $g \ra \infty$ is not of Tate type.

Universiteit Utrecht, Mathematisch Instituut, Postbus 80.010, 3508 TA Utrecht,
Nederland; e-mail: pikaart@math.ruu.nl

\end{document}